
\def\LaTeX{%
  \let\Begin\begin
  \let\End\end
  \let\salta\relax
  \let\finqui\relax
  \let\futuro\relax}

\def\UK{\def\our{our}\let\sz s}
\def\USA{\def\our{or}\let\sz z}



\LaTeX

\USA


\salta

\documentclass[twoside,12pt]{article}
\setlength{\textheight}{23.5cm}
\setlength{\textwidth}{16cm}
\setlength{\oddsidemargin}{2mm}
\setlength{\evensidemargin}{2mm}
\setlength{\topmargin}{-15mm}
\parskip2mm


\usepackage{color}
\usepackage{amsmath}
\usepackage{amsthm}
\usepackage{amssymb}
\usepackage[mathcal]{euscript}

 


\definecolor{dbrown}{rgb}{0.5,0.25,0}



\bibliographystyle{plain}


%

\finqui

\def\Beq{\Begin{equation}}
\def\Eeq{\End{equation}}
\def\Bsist{\Begin{eqnarray}}
\def\Esist{\End{eqnarray}}

\def\Bthm{\Begin{theorem}}
\def\Ethm{\End{theorem}}
\def\Blem{\Begin{lemma}}
\def\Elem{\End{lemma}}
\def\Bprop{\Begin{proposition}}
\def\Eprop{\End{proposition}}

\def\Brem{\Begin{remark}\rm}
\def\Erem{\End{remark}}

\let\non\nonumber



\newcommand\QED{\hfill $\square$}


\def\step #1 \par{\medskip\noindent{\bf #1.}\quad}


\def\holder{H\"older}
\def\aand{\quad\hbox{and}\quad}
\def\loti{long-time}
\def\Loti{Long-time}

\def\lhs{left-hand side}
\def\rhs{right-hand side}

\def\wk{well-known}
\def\wepo{well-posed}
\def\Wepo{Well-posed}


\def\characteriz{characteri\sz}
\def\generaliz{generali\sz}

\def\organiz{organi\sz}

\def\bhv{behavi\our}


\def\multibold #1{\def\arg{#1}%
  \ifx\arg\pto \let\next\relax
  \else
  \def\next{\expandafter
    \def\csname #1#1#1\endcsname{{\bf #1}}%
    \multibold}%
  \fi \next}

\def\pto{.}

\def\multical #1{\def\arg{#1}%
  \ifx\arg\pto \let\next\relax
  \else
  \def\next{\expandafter
    \def\csname cal#1\endcsname{{\cal #1}}%
    \multical}%
  \fi \next}


\def\multimathop #1 {\def\arg{#1}%
  \ifx\arg\pto \let\next\relax
  \else
  \def\next{\expandafter
    \def\csname #1\endcsname{\mathop{\rm #1}\nolimits}%
    \multimathop}%
  \fi \next}

\multibold
qwertyuiopasdfghjklzxcvbnmQWERTYUIOPASDFGHJKLZXCVBNM.

\multical
QWERTYUIOPASDFGHJKLZXCVBNM.

\multimathop
dist div dom meas sign supp .


\def\accorpa #1#2{\eqref{#1}--\eqref{#2}}
\def\Accorpa #1#2 #3 {\gdef #1{\eqref{#2}--\eqref{#3}}%
  \wlog{}\wlog{\string #1 -> #2 - #3}\wlog{}}


\def\graffe #1{\mathopen\{#1\mathclose\}}

\def\<#1>{\mathopen\langle #1\mathclose\rangle}
\def\norma #1{\mathopen \| #1\mathclose \|}

\def\iot {\int_0^t}
\def\ioT {\int_0^T}
\def\ioi {\int_0^\infty}
\def\iO{\int_\Omega}
\def\intQt{\iot\!\!\iO}
\def\intQ{\ioT\!\!\iO}
\def\intQi{\ioi\!\!\iO}

\def\dt{\partial_t}
\def\dn{\partial_\nu}

\def\cpto{\,\cdot\,}

\def\checkmmode #1{\relax\ifmmode\hbox{#1}\else{#1}\fi}
\def\aeO{\checkmmode{a.e.\ in~$\Omega$}}
\def\aeQ{\checkmmode{a.e.\ in~$Q_T$}}

\def\aaQ{\checkmmode{for a.a.~$(x,t)\in Q_T$}}
\def\aat{\checkmmode{for a.a.~$t\in(0,T)$}}
\def\Aat{\checkmmode{for a.a.~$t>0$}}


\def\erre{{\mathbb{R}}}




\def\genspazio #1#2#3#4#5{#1^{#2}(#5,#4;#3)}
\def\spazio #1#2#3{\genspazio {#1}{#2}{#3}T0}

\def\L {\spazio L}
\def\H {\spazio H}
\def\W {\spazio W}

\def\C #1#2{C^{#1}([0,T];#2)}


\def\Lx #1{L^{#1}(\Omega)}
\def\Hx #1{H^{#1}(\Omega)}
\def\Wx #1{W^{#1}(\Omega)}

\def\Luno{\Lx 1}
\def\Ldue{\Lx 2}
\def\Linfty{\Lx\infty}
\def\Lq{\Lx4}
\def\Huno{\Hx 1}
\def\Hdue{\Hx 2}


\def\LQ #1{L^{#1}(Q_T)}


\let\theta\vartheta
\let\eps\varepsilon
\let\phi\varphi

\let\TeXchi\chi                         
\newbox\chibox
\setbox0 \hbox{\mathsurround0pt $\TeXchi$}
\setbox\chibox \hbox{\raise\dp0 \box 0 }
\def\chi{\copy\chibox}


\def\muz{\mu_0}
\def\rhoz{\rho_0}

\def\rhomin{\rho_*}

\def\normaV #1{\norma{#1}_V}
\def\normaH #1{\norma{#1}_H}

\def\mun{\mu_n}
\def\rhon{\rho_n}

\def\un{u_n}
\def\ui{u_\infty}

\def\rhomrm{(\rho-\rmin)^-}
\def\rhomrmeps{(\rhoeps-\rmin)^-}
\def\rmin{r_*}

\def\mui{\mu_\infty}
\def\rhoi{\rho_\infty}

\def\mus{\mu_s}
\def\rhos{\rho_s}
\def\muo{\mu_\omega}
\def\rhoo{\rho_\omega}
\def\mueps{\mu_\eps}
\def\rhoeps{\rho_\eps}
\def\ueps{u_\eps}

\def\Vp{V^*}

\def\neto{\mathrel{\scriptstyle\nearrow}}
\def\seto{\mathrel{\scriptstyle\searrow}}

\DeclareMathAlphabet{\mathbf}{OT1}{cmr}{bx}{it}

\newcommand{\hb}{\mathbf{h}}
\newcommand{\nb}{\mathbf{n}}
\newcommand{\Mb}{\mathbf{M}}
\newcommand{\Hb}{\mathbf{H}}
 \font\mba=cmmib10 scaled
\magstephalf

\def\csib{\hbox{\mba {\char 24}}}
\Begin{document}


\title{{\bf An asymptotic analysis for a nonstandard Cahn-Hilliard system
with viscosity}\footnote{{\bf Acknowledgments.}\quad\rm
The authors gratefully acknowledge the financial
support of the MIUR-PRIN Grant 2008ZKHAHN \emph{``Phase transitions, hysteresis 
and multiscaling''}, the EU Marie Curie Research Training Network
MULTIMAT \emph{``Multi-scale Modeling and Characterization
for Phase Transformations in Advanced Materials''},
the DFG Research Center {\sc Matheon} in Berlin, 
and the IMATI of CNR in Pavia.}}
\author{}
\date{}
\maketitle
\begin{center}
\vskip-2.2cm
{\large\bf Pierluigi Colli$^{(1)}$}\\
{\normalsize e-mail: {\tt pierluigi.colli@unipv.it}}\\[.25cm]
{\large\bf Gianni Gilardi$^{(1)}$}\\
{\normalsize e-mail: {\tt gianni.gilardi@unipv.it}}\\[.25cm]
{\large\bf Paolo Podio-Guidugli$^{(2)}$}\\
{\normalsize e-mail: {\tt ppg@uniroma2.it}}\\[.25cm]
{\large\bf J\"urgen Sprekels$^{(3)}$}\\
{\normalsize e-mail: {\tt sprekels@wias-berlin.de}}\\[.45cm]
$^{(1)}$
{\small Dipartimento di Matematica ``F. Casorati'', Universit\`a di Pavia}\\
{\small via Ferrata 1, 27100 Pavia, Italy}\\[.2cm]
$^{(2)}$
{\small Dipartimento di Ingegneria Civile, Universit\`a di Roma ``Tor Vergata''}\\
{\small via del Politecnico 1, 00133 Roma, Italy}\\[.2cm]
$^{(3)}$
{\small Weierstra\ss-Institut f\"ur Angewandte Analysis und Stochastik}\\
{\small Mohrenstra\ss e\ 39, 10117 Berlin, Germany}\\[.8cm]
\end{center}

\Begin{abstract}
This paper is concerned with a diffusion model of phase-field type, consisting 
of a {parabolic} system of two partial differential equations{,} interpreted as balances
of microforces and microenergy{, for two unknowns: the problem's order parameter $\rho$} 
and the chemical potential {$\mu$; each equation  includes a viscosity term -- respectively, $\varepsilon \,\partial_t\mu$ and $\delta\,\partial_t\rho$ -- with $\varepsilon$ and $\delta$ two positive parameters; the field equations are complemented by Neumann homogeneous boundary conditions and suitable initial conditions.  In a recent paper \cite{CGPS3}, we proved that this problem is  \wepo\ and investigated the \loti\ \bhv\ of its $(\varepsilon,\delta)-$solutions. Here we  discuss the asymptotic limit of the system as $\eps$ 
tends to $0$. We prove convergence of  $(\varepsilon,\delta)-$solutions to the corresponding solutions for the case $\eps =0$, whose long-time behavior we characterize; in the proofs, we employ compactness and monotonicity arguments.} 
\End{abstract}

\salta

\pagestyle{myheadings}
\newcommand\testopari{\sc Colli \ --- \ Gilardi \ --- \ Podio-Guidugli \ --- \ Sprekels}
\newcommand\testodispari{\sc Asymptotic analysis for a nonstandard diffusion model}
\markboth{\testodispari}{\testopari}

\finqui


\section{Introduction}
\label{Intro}
\setcounter{equation}{0}

{The system we study was proposed for mathematical investigation in \cite{CGPS3}; as to modeling issues, its most directly relevant antecedents are two papers by Fried \& Gurtin \cite{FG} and Gurtin \cite{Gurtin}, and a paper by one of us \cite{Podio}.

\step A nonstandard phase-field evolution problem

The initial/boundary-value problem we dealt with in \cite{CGPS3} consists in finding two \emph{phase fields},  the {\it chemical potential} $\mu$ and the \emph{order parameter} $\rho$,  such that
\Bsist
  & \eps\, \dt\mu + 2\rho \,\dt\mu + \mu \, \dt\rho - \Delta\mu = 0
  & \quad \hbox{in $\Omega \times (0,+\infty)$,}
  \label{Iprima}
  \\
  & \delta\, \dt\rho - \Delta\rho + f'(\rho) = \mu 
  & \quad \hbox{in $\Omega \times (0,+\infty)$,}
  \label{Iseconda}
  \\
  & \dn\mu = \dn\rho = 0
  & \quad \hbox{on $\Gamma \times (0,+\infty)$,}
  \label{Ibc}
  \\
  & \mu(\cpto,0) = \muz
  \aand
  \rho(\cpto,0) = \rhoz
  & \quad \hbox{in $\Omega$,}
  \label{Icauchy}
\Esist
\Accorpa\Ipbl Iprima Icauchy
where $\Omega$ denotes a bounded domain of $ \erre^3$ with (sufficiently) smooth boundary $\Gamma$, and $f'$ stands for 
the derivative of a double-well potential $f$. 
 This nonstandard   phase-field model
can be regarded as a variant of the classic Cahn-Hilliard system for diffusion-driven phase segregation by atom rearrangement:
\Beq\label{CH}
 \dt\rho - \kappa \Delta \mu =0 \ , \qquad  \mu= - \Delta\rho + 
f'(\rho)
 {.}  
 \Eeq
{Note, in \eqref{Iprima}, the} 
unpleasant nonlinear terms involving time derivatives, and the fact that we have taken the \emph{mobility} coefficient $\kappa >0 $ equal to $1$. Moreover, recall that  
equations \eqref{CH} are customarily combined so as to obtain the well-known 
\emph{Cahn-Hilliard equation}:
\Beq\label{CHe}
\dt\rho = \kappa \Delta (- \Delta\rho + f'(\rho)).
\Eeq
%

\step Fried \& Gurtin's generalization of Cahn-Hilliard equation

In \cite{FG,Gurtin} a broad generalization of \eqref{CHe} was arrived at, with three measures:
\begin{enumerate}
 \item[(i)]
by regarding the second of \eqref{CH} as a \textit{balance of microforces:}
\Beq\label{balance}
\div\csib+\pi+\gamma=0{,}
\Eeq
where the distance microforce per unit volume 
is split into an internal part $\pi$ and an 
external part $\gamma$, and the contact microforce per unit area of a surface oriented by its normal $\nb$ is measured by $\csib\cdot\nb$ in terms of the \emph{microstress} vector  $\csib$;\footnote{In \cite{Fremond}, the microforce balance
is stated under form of a principle of virtual powers 
for microscopic motions.}

\item[(ii)] 
by interpreting the first equation of \eqref{CH} as a 
\textit{balance law for the order parameter}:
\Beq\label{balorpam}
\partial_t\rho = - \div \hb + \sigma{,}
\Eeq
where the pair $(\hb ,  \sigma)$ is the \textit{inflow} of $\rho$; 

\item[(iii)] 
by requiring that the constitutive choices for $\pi,\csib, \hb$, and the \emph{free energy density} $\psi$, {be} consistent in the sense of Coleman and Noll \cite{CN} with a postulated ``dissipation inequality that accomodates diffusion'':
\begin{equation}\label{dissipation}
\partial_t\psi +(\pi-\mu)\partial_t\rho-\csib\cdot\nabla(\partial_t\rho)+\hb\cdot\nabla\mu\leq 0
\end{equation}
(cf. equation (3.6) in \cite{Gurtin}). 
\end{enumerate}
In \cite{Gurtin}, the following set of constitutive prescriptions was shown to be consistent with (iii):
\Beq\label{costi} \left\{
\begin{array}{c}
\psi = \widehat\psi(\rho,\nabla\rho),  \displaystyle \\[0.1cm]
\widehat\pi(\rho,\nabla\rho,\mu)=\mu- \partial_\rho \widehat\psi(\rho,\nabla\rho), \displaystyle \\[0.1cm] 
\widehat\csib(\rho,\nabla\rho)=\partial_{\nabla\rho} \widehat\psi(\rho,\nabla\rho) \displaystyle
\end{array}
\right\}
\Eeq 
 together with 
 \Beq\label{acca}
\hb = - \Mb\nabla \mu , \quad \hbox{with } \ \Mb=\widehat\Mb(\rho,\nabla\rho,\mu, \nabla\mu),
\Eeq
provided the tensor-valued \emph{mobility mapping} $\widehat\Mb$ satisfies the  \emph{residual dissipation
{inequality}}
\Beq
\nabla \mu\cdot \widehat\Mb(\rho,\nabla\rho,\mu, \nabla\mu) \nabla\mu \geq 0 . \non
\Eeq
With the use of \eqref{balance}, \eqref{balorpam}, \eqref{costi}, and $\eqref{acca}_1$, a general equation for diffusive phase segregation processes is arrived at: 
%
\[
\dt\rho = \div\left(\Mb\nabla\left(\partial_\rho \widehat\psi(\rho,\nabla\rho)-\div\big(\partial_{\nabla\rho} \widehat\psi(\rho,\nabla\rho)\big)-\gamma\right)\right)+\sigma {;}
\]
in particular, the Cahn-Hilliard equation \eqref{CHe} is obtained by {taking} the external distance microforce $\gamma$ and the order-parameter source term $\sigma$ identically null, and by choosing
\Bsist
&&\widehat\psi(\rho,\nabla\rho)= f(\rho)+\frac{1}{2}|\nabla\rho|^2,\qquad \Mb=\kappa \mathbf{1}.
\label{constitutive}
\Esist

\step An alternative generalization of Cahn-Hilliard equation

In \cite{Podio}, a major modification of Fried \& Gurtin's approach to phase-segregation modeling was proposed. While the crucial step (i) was retained, both the order-parameter balance \eqref{balorpam} and the dissipation inequality \eqref{dissipation} were} dropped and replaced, respectively, by
 the \emph{microenergy balance}
\Beq\label{energy}
\partial_t\varepsilon=e+w,\quad e:=-\div{\overline \hb}+{\overline \sigma},\quad w:=-\pi\,\partial_t\rho+\csib\cdot\nabla(\partial_t\rho){,}
\Eeq
and the \emph{microentropy imbalance}
\Beq\label{entropy}
\partial_t\eta\geq -\div\hb+\sigma,\quad \hb:=\mu{\overline \hb},\quad \sigma:=\mu\,{\overline \sigma}.
\Eeq
{A further new feature was} that the \emph{microentropy inflow} $(\hb,\sigma)$ was deemed proportional to the \emph{microenergy inflow} $({\overline \hb},{\overline\sigma})$ through the \emph{chemical potential} $\mu$, a {positive} field; consistently, the free energy was defined~to~be
\Beq\label{freeenergy}
\psi:=\varepsilon-\mu^{-1}\eta,
\Eeq
with {the} chemical potential playing the same role as 
{the} \emph{coldness} in the deduction of the heat equation.\footnote{As much as absolute temperature is a macroscopic measure of microscopic 
\emph{agitation}, its inverse - the coldness - measures microscopic \emph{quiet}; likewise, as argued in 
\cite{Podio}, {the} chemical potential can be seen as a macroscopic measure of microscopic \emph{organization}.}  

Combination of (\ref{energy})-(\ref{freeenergy}) gives:
\Beq\label{reduced}
\partial_t\psi\leq -\eta_{}\dt (\mu^{-1})+\mu^{-1}\,{\overline \hb}\cdot\nabla\mu-\pi\,\partial_t\rho+\csib\cdot\nabla(\partial_t\rho),
\Eeq
an inequality that replaces (\ref{dissipation}) in restricting  \emph{\`a la} Coleman-Noll {the possible} constitutive choices. 
On taking all of the constitutive mappings delivering $\pi,\csib,\eta$, and ${\overline \hb}$, 
depend{ent in principle} on  $\rho,\nabla\rho,\mu,\nabla\mu$, and on choosing
\Beq\label{constitutives}
\psi=\widehat\psi(\rho,\nabla\rho,\mu)=-\mu\,\rho+f(\rho)+\frac{1}{2}|\nabla\rho|^2,
\Eeq
compatibility with (\ref{reduced}) implies that we must have:
\Beq \label{cn}
\left\{
\begin{array}{c}
\widehat\pi(\rho,\nabla\rho,\mu)=-{\partial_{{\rho}} \widehat\psi(\rho,\nabla\rho,\mu)}~\displaystyle{=\mu-f'(\rho),}\\[0.2cm]
\widehat\csib(\rho,\nabla\rho,\mu)={\partial_{{\nabla\rho}} \widehat\psi(\rho,\nabla\rho,\mu)}=\nabla\rho, \displaystyle\\[0.2cm]
\widehat\eta(\rho,\nabla\rho,\mu)=\mu^2 \partial_{{\mu}} \widehat\psi(\rho,\nabla\rho,\mu)\displaystyle{=-\mu^2\rho}  
\end{array}
\right\}
\Eeq
together with
\Beq
\widehat{\overline \hb}(\rho,\nabla\rho,\mu,\nabla\mu) =  \widehat\Hb(\rho,\nabla\rho,\mu, \nabla\mu)\nabla \mu , \quad
\nabla \mu\cdot \widehat\Hb(\rho,\nabla\rho,\mu, \nabla\mu) \nabla\mu \geq 0 . \non
\Eeq
If {we now choose} for $\widehat\Hb$ the simplest expression $\Hb=\kappa \mathbf{1} $, implying a constant {and isotropic} mobility, and {if we once again} assume that the external distance microforce $\gamma$ and the source $\overline \sigma$ are null, {then}, with the use of (\ref{cn}) and \eqref{freeenergy}, the microforce balance (\ref{balance}) and the energy balance (\ref{energy}) become, respectively, 
\Beq\label{a}
 \Delta\rho+\mu- f'(\rho)  = 0
\Eeq
and
\Beq
 2\rho \,\dt\mu + \mu \, \dt\rho - \kappa \Delta\mu = 0,   \label{secondanew}
\Eeq
a nonlinear system for the unknowns $\rho$ and $\mu$. 

\step {Two-parameter regularization}

{Compare now the systems \eqref{a}-\eqref{secondanew}
and \eqref{CH}: \eqref{a} and $\eqref{CH}_2$ are one and the same `static' relations  between $\mu $ and $\rho$, whereas} \eqref{secondanew}
is rather different from $\eqref{CH}_1$, for more than one reason: 
\begin{enumerate}

\item[(R1)]  
$\eqref{CH}_1$ is linear,
\eqref{secondanew} is not;
  
\item [(R2)] the time derivatives of $\rho$ and $\mu$ 
are both present in  \eqref{secondanew}; 

\item[(R3)] 
in front of both $\dt\mu $ and 
$\dt\rho$ there are {\it nonconstant factors} that {\it should remain nonnegative} during the evolution. 

\end{enumerate} 

Thus, the system \eqref{a}-\eqref{secondanew} deserves a careful 
{analysis}. We must confess that we {boldly}
 attacked this problem as is, 
prompted to optimism by the successful outcome of a previous joint 
research effort {\cite{CGPS, CGPS2} devoted to tackling the system of Allen-Cahn type one arrives at via the 
approach in \cite{Podio} for no-diffusion phase-segregation processes.}
Unfortunately, {the evolution problem ruled by} \accorpa{a}{secondanew} turned out to be too 
difficult for us. Therefore, we decided to study its regularized 
version \Ipbl, {which we} obtained by introducing the extra terms $\eps\, 
\partial_t\mu $ in \eqref{secondanew} and $\delta\,\partial_t\rho$ 
in \eqref{a}, for small positive coefficients $\eps$ and 
$\delta$. Motivations for the introduction of such terms are proposed 
and discussed in \cite{CGPS3}; {interestingly, while the second can be interpreted as a dissipative part of the distance microforce, so far we have not been able to find a convincing physical interpretation for the first. But, our present study demonstrates -- so we believe -- that it can legitimately be regarded as an efficient mathematical device.}

\step Limit as the first parameter tends to 0

In \cite{CGPS3}, by assuming (as we did in  \cite{CGPS, CGPS2})
that $f'$ is the sum of a strictly increasing $C^1$ function 
$f'_1$ with domain $(0,1)$ {that is singular at the} endpoints, 
and of a smooth bounded perturbation $f'_2$
(to allow for a double- or multi-well potential $f$), 
we first proved existence of a strong solution $(\mu, \rho)$ 
to \Ipbl\ satisfying  $\mu \geq 0$ and 
$0< \rho <1 $ almost everywhere in $\Omega \times (0,+\infty)$ 
(of course, {we stipulated that the initial data} meet same requirements in $\Omega$).   
{Then,} under {some additional technical} assumptions, we showed that the 
component $\mu$ is bounded, and {so is} $f'(\rho) $; as a consequence, $\rho$ 
{stays away from the threshold values $0$ and $1$}.
These boundedness properties are very useful in proving uniqueness of solutions.

In some sense, passing to the limit as the regularizing parameters tend to zero is
the challenging final aim of our {research} project. 
For the moment {being},
we are able to deal with $\eps$ and to deduce, by {a rather delicate} asymptotic analysis, an existence theorem for the {limit} problem. {Precisely,} we let $\eps$ tend to zero and show that any weak or weak star limit
of {any} subsequence of solutions $(\mu_\eps, \rho_\eps)$ 
to \Ipbl\ yields a solution $(\mu, \rho)$
to the resulting limit problem, {which is} obtained by putting $\eps=0$ in \Ipbl\ and 
rewriting the corresponding first equation in the form
\Beq
  2 \dt(\mu\rho) - \Delta\mu = \mu \, \dt\rho
  \quad \hbox{rather than as} \quad
  2 \rho \, \dt\mu + \mu \, \dt\rho - \Delta\mu = 0 .
  \label{Ilimprima}
\Eeq
{This we do because} {it is not clear} {to us} {from the structure of the system} 
 whether {a suitably regular representation} for $\dt\mu$ {could be recovered} in the limit,
while {we are able to} show that the time derivative $\dt(\mu\rho)$ actually exists, 
at least in some \generaliz ed sense.\footnote{In fact, \eqref{Iprima} has the equivalent formulation: 
\[
  \dt (\eps\mu + 2\mu\rho) - \Delta\mu = \mu \dt {\rho,} 
\]
which singles out the time derivative of the auxiliary variable $\eps\mu + 2\mu\rho$ for $\eps>0$.}


Here,  {as we did in \cite{CGPS3},} we also deal with the long-time behavior of the system. 
We prove that each element $(\mu_\omega, \rho_\omega)$ of the $\omega$-limit set for a certain trajectory
is a steady state solution {to} \Ipbl; therefore, in particular, $\mu_\omega$ \emph{is a constant} (cf. \eqref{Ilimprima} and \eqref{Ibc}).\footnote{Note {that} the steady state problem associated with both cases $\eps>0$ and $\eps =0$ is the same.}  

An outline of our paper is the following: in Section~\ref{MainResults}, we {carefully state} assumptions and 
results; Section~\ref{ConvProof} is devoted to the proof of the convergence theorem, so as to deduce the 
existence of solutions {to} the limit problem; finally, in Section~\ref{OmegaProof}, we develop our argument 
for the characterization of the $\omega$-limit.


\section{{Assumptions and} main results}
\label{MainResults}
\setcounter{equation}{0}

First of all,
we assume $\Omega$ to be a bounded connected open set in $\erre^3$
with smooth boundary~$\Gamma$
and set, for convenience,
\Beq
  V := \Huno,
  \quad H := \Ldue ,
  \aand
  W := \graffe{v\in\Hdue:\ \dn v = 0}.
  \label{defspazi}
\Eeq
We endow 
 {these} spaces with their standard norms,
for which we use self-explaining notation like~$\normaV\cpto$.
However, we write $\normaH\cpto$ for the norm in any power of $H$ as well.
The symbol $\<\cpto,\cpto>$ denotes the duality product 
between~$\Vp$, the dual space of~$V$, and~$V$ itself.
Since $\Omega$ is bounded and smooth, the embeddings $W\subset V\subset H$ are compact.
Moreover, since $V$ is dense in~$H$,
we can identify $H$ with a subspace of~$\Vp$ 
in the usual~way,
i.e., in order that 
 {$\<u,v>=(u,v)_H$, where $(\,\cdot\,,\,\cdot)_H$
denotes the inner product in $H$, holds for every $u\in H$ and $v\in V$.}
Then, also the embedding $H\subset\Vp$ is compact.

As in~\cite{CGPS3}, we assume~that
\Bsist
  && f = f_1 + f_2,
  \quad \hbox{where} \quad
  \hbox{$f_1,f_2:(0,1) \to \erre$ are functions satisfying:}
  \label{hpf}
  \\
  && \hbox{$f_1$ is $C^1$ and convex}, \quad
  \hbox{$f_2$ is $C^2$},
  \aand
  \hbox{$f_2''$ is bounded}{,}
  \label{hpfi}
  \\
  && \lim_{r\searrow0} f_1'(r) = - \infty
  \aand
  \lim_{r\nearrow0} f_1'(r) = + \infty .
  \label{hpfp}
\Esist
\Accorpa\Hpf hpf hpfp
As to initial data, we start with the assumptions
\Bsist
  && \muz \in V
  \aand \muz \geq 0 \quad \aeO;
  \label{hpmuz}
  \\
  && \rhoz \in W, \quad
  0 < \rhoz < 1 \quad \hbox{in $\Omega$};
  \aand
  f'(\rhoz) \in H .
  \label{hprhoz}
\Esist
\Accorpa\Hpdati hpmuz hprhoz
{The reader is referred to} the forthcoming Remark~\futuro\ref{Hpdeboli}
for weaker conditions.

 {Since we aim to let}
$\eps$ tend to zero,
we stress the dependence of the solution
found in~\cite{CGPS3} on the parameter~$\eps$.
In that paper, for any fixed $T>0$, the following a~priori regularity is required:
\Bsist
  && \mueps \in \H1H \cap \L2W,
  \label{regmueps}
  \\
  && \rhoeps \in \W{1,\infty}H \cap \H1V \cap \L\infty W,
  \label{regrhoeps}
  \\
  && \mueps \geq 0 \quad \aeQ,
  \label{muepspos}
  \\
  && 0 < \rhoeps <1 \quad \aeQ
  \aand
  f'(\rhoeps) \in \L\infty H ,
  \label{regfprhoeps}
\Esist
\Accorpa\Regsoluzeps regmueps regfprhoeps
where $Q_T:=\Omega\times(0,T)$.
We note that homogeneous Neumann boundary conditions follow from
\accorpa{regmueps}{regrhoeps}, in view of the definition of~$W$
(see~\eqref{defspazi}).
Then, the $\eps$-problem is:
\Bsist
  & {\Big[}(\eps + 2\rhoeps) \dt\mueps + \mueps \, \dt\rhoeps - \Delta\mueps = 0{\quad \hbox{or}\Big]}
  &
  \non
  \\
  & 
  \dt (\eps\mueps + 2\mueps\rhoeps) - \Delta\mueps = \mueps \dt\rhoeps
  & \quad \aeQ,
  \label{primaeps}
  \\
  & \delta \dt\rhoeps - \Delta\rhoeps + f'(\rhoeps) = \mueps
  & \quad \aeQ,
  \label{secondaeps}
  \\
  & \mueps(0) = \muz
  \aand
  \rhoeps(0) = \rhoz
  & \quad \aeO .
  \label{cauchyeps}
\Esist
\Accorpa\Pbleps primaeps cauchyeps

We recall the existence result of~\cite{CGPS3}.

\Bthm
\label{Existence}
Let $T\in(0,+\infty)$, and assume that \Hpf\ and \Hpdati\ are satisfied.
Then, there exists a pair $(\mueps,\rhoeps)$ 
satisfying \Regsoluzeps\ and solving problem~\Pbleps.
\Ethm

\noindent{As to uniqueness, the result in {\cite[Thm.~2.2]{CGPS3}}  holds for solutions that, in addition 
to \Regsoluzeps, have certain properties that, in turn, are guaranteed}
whenever the {initial} data fulfil the following conditions, {additional to \eqref{hpmuz} and \eqref{hprhoz}}:
$$\muz\in\Linfty, \quad\inf\rhoz>0, \quad\textrm{and}\quad\sup\rhoz<1$$
(see \cite[Thm~2.3]{CGPS3}).
Within such a framework,
since $T>0$ is arbitrary,
the existence of a unique solution $(\mueps,\rhoeps)$ 
defined for every positive time was ensured,
and its \loti\ \bhv\ was studied.

 {Here,} {our first concern is}
to construct a global solution $(\mueps,\rhoeps)$
to problem \Pbleps\ that satisfies \Regsoluzeps\ for every finite~$T$,
without assuming the just mentioned stronger conditions on the initial data.
We 
 {cannot}
 ensure uniqueness, of course. 
 {The corresponding result reads:}

\Bprop
\label{Globalexistence}
Assume that \Hpf\ and \Hpdati are fulfilled.
Then, there exists a pair $(\mueps,\rhoeps):[0,+\infty)\to W\times W$ 
satisfying \Regsoluzeps\ and solving
problem~\Pbleps\ for every $T\in(0,+\infty)$.
\Eprop

 {Starting from any family $\{(\mueps,\rhoeps)\}_{\eps >0}$ of 
solutions of this type, we then let $\eps$ tend to zero. To do this, we need to assume that}
\Beq
  \inf\rhoz > 0 ,
  \label{hpinfrhoz}
\Eeq
in addition to \Hpdati.
 {Under this assumption, we show that}
$\rhoeps$ is bounded away from zero and that
$(\mueps,\rhoeps)$ tends to some pair $(\mu,\rho)$ as $\eps\seto0$
in a suitable topology, at least for a subsequence.
Moreover, we determine the limit problem solved by~$(\mu,\rho)$.
The a~priori regularity we require for $(\mu,\rho)$ 
on every finite time interval is the following:
\Bsist
  && \mu \in \L\infty H\cap \L2V,
  \label{regmu}
  \\
  && \rho \in \H1H \cap \L\infty V \cap \L2W,
  \label{regrho}
  \\
  && \mu \geq 0 
  \quad \aeQ
  \aand
  \inf\rho > 0,
  \label{murhopos}
  \\
  && 
  \rho <1
  \quad \aeQ
  \aand
  f'(\rho) \in \L2H,
  \label{regfprho}
  \\
  && \mu\rho \in \W{1,p}\Vp 
  \quad \hbox{for some $p\in(1,+\infty)$};
  \label{regprod}
\Esist
\Accorpa\Regsoluz regmu regprod
the 
 {corresponding limit problem is:}
\Bsist
  && 2 \< \dt (\mu\rho)(t) , v >
  + \iO \nabla\mu(t) \cdot \nabla v
  = \iO \mu(t) \, \dt\rho(t) \, v
  \non
  \\
  && \quad \hbox{for every $v\in V$ and \aat},
  \label{prima}
  \\
  && \delta \dt\rho - \Delta\rho + f'(\rho) = \mu
  \quad \aeQ,
  \vrule depth 5pt width0pt
  \label{seconda}
  \\
  && (\mu\rho)(0) = \muz\rhoz
  \aand
  \rho(0) = \rhoz
  \quad \aeO .
  \label{cauchy}
\Esist
\Accorpa\Pbl prima cauchy

\Brem
\label{Commentoprima}
{The last integral in \eqref{prima} makes sense because}
$V\subset\Lx3$ and
$\mu\dt\rho$ belongs at least to $\L1{\Lx{3/2}}$,
as a consequence of \accorpa{regmu}{regrho}.
{Note that
\eqref{prima} incorporates} the homogeneous Neumann boundary condition for $\mu$
in a \generaliz ed sense. {Moreover, note}
that \eqref{regprod} implies that
$\mu\rho$ is a continuous $\Vp$-valued function,
so that the first equality in \eqref{cauchy} has a precise meaning.
On the contrary, no continuity for $\mu$ is ensured at the moment.
\Erem

Here is our convergence result.

\Bthm
\label{Convergence}
Assume that \Hpf, \Hpdati, and~\eqref{hpinfrhoz} are satisfied,
and let $\graffe{(\mueps,\rhoeps)}_{\eps\in(0,1)}$
be a family of solutions to problem \Pbleps\ satisfying \Regsoluzeps.
Then, there exists a pair $(\mu,\rho)$, satisfying \Regsoluz\ and solving problem \Pbl\
for every $T\in(0,+\infty)$,
such that $(\mueps,\rhoeps)$ converges to $(\mu,\rho)$
in a suitable topology,\footnote{{to}
 {be specified in the course of}
the proof given in the next section.}
at least for a subsequence $\eps_n\seto0$.
\Ethm

\Brem
\label{Hpdeboli}
The assumptions \Hpdati\ are strong. In fact, while they are needed 
for Theorem~\ref{Existence}
and Proposition~\ref{Globalexistence},
some of them will not play any role in the 
 {following,
as it can be seen}
by looking at our a priori estimates.
For instance, the last 
 {condition in}
\eqref{hprhoz} will not be important,
since just $f(\rhoz)\in\Luno$ will be used.
{Accordingly}, one can prove a result similar to Theorem~\ref{Convergence},
{but involving} $\eps$-approximations of less regular initial data that{, this notwithstanding, satisfy} \Hpdati.
{Precisely, suppose we assume that} 
\Beq
  \muz \in H , \quad
  \muz \geq 0 \quad \aeO , \quad
  \rhoz \in V , \quad
  f(\rhoz) \in \Luno, 
  \aand
  \inf\rhoz > 0 .
  \label{hpdeboli}
\Eeq
{Then, it would be} possible to construct $\eps$-approximations
 {$(\mu_{0\eps},\rho_{0\eps})$ of such 
initial data $(\muz,\rhoz)$ that satisfy \Hpdati\, and whose norms of type \eqref{hpdeboli} 
remain bounded as $\eps\seto0\,$:}
{e.g., as to} $\rho_{0\eps}$,
one {could} take the solution to the elliptic equation
\Beq
  \frac {\rho_{0\eps} - \rhoz} \eps
  - \Delta\rho_{0\eps} + f_1'(\rho_{0\eps})
  = 0
  \quad \aeO,
  \non
\Eeq
{supplemented by} homogeneous Neumann boundary conditions.
\Erem

\Brem
\label{Esdebole} 
Theorem~\ref{Convergence} and the previous Remark offer us the 
possibility of defining and obtaining a weaker solution to problem \Pbleps\  
(that is, also for the case $\eps >0$), {if one writes}
{equation \eqref{Iprima} in the form \eqref{primaeps}.}
To see that this solution is weaker than the one provided by Theorem~\ref{Existence}, it suffices to compare \accorpa{regmu}{regrho} with 
\accorpa{regmueps}{regrhoeps}. On the other hand, we can just assume \eqref{hpdeboli} 
and point out that in this approach one should consider \eqref{regprod}, 
\eqref{prima}, \eqref{cauchy} with $\mu\rho$ replaced by $(\eps/2)\mu + \mu\rho $.
\Erem

Our  {final} aim is to study the \loti\ \bhv\
of any solution constructed {according to} Theorem~\ref{Convergence}.
To this end, we introduce the $\omega$-limit of the trajectory
in a proper topology,
and prove that every element of it coincides with a steady state.
We~set:
\Bsist
  & \omega(\mu,\rho)
  & = \bigl\{
    (\muo,\rhoo) \in H \times V:\ \bigl( \mu(t_n),\rho(t_n) \bigr) \to (\muo,\rhoo)\ 
  \non
  \\
  && \qquad
    \hbox{weakly in $H\times V$ for some sequence $t_n\neto+\infty$}
  \bigr\} .
  \label{defomegalim}
\Esist
The above definition has a precise meaning, {because the pointwise values of the $(H\times V)$-valued function $(\mu,\rho)$ are well defined thanks to the continuity properties stated in our next result.}

\Bthm
\label{Omegalimit}
Assume that \Hpf, \Hpdati\, and~\eqref{hpinfrhoz} are satisfied,
and let $(\mu,\rho):[0,+\infty)\to H\times V$ be 
given by Theorem~\ref{Convergence}.
Then, $(\mu,\rho)$ is bounded, 
and its components $\mu$ and $\rho$ are weakly and strongly continuous, respectively.
In particular, the $\omega$-limit \eqref{defomegalim}
is nonempty.
Moreover, every element of $\omega(\mu,\rho)$
coincides with a pair $(\mus,\rhos)$
 {such that}
\Bsist
  && \hbox{$\mus$~is a nonnegative constant,}
  \non
  \\
  && \rhos \in W , \quad
  0 < \rhos < 1 , \quad
  f'(\rhos) \in H ,
  \aand
  -\Delta\rhos + f'(\rhos) = \mus
  \quad \aeO, \qquad\quad 
  \label{secondas}
\Esist
i.e., {it coincides} with a steady state.
\Ethm

We stress that the above result does not 
{necessarily {hold for} all possible solutions.}
Indeed, it 
 {only deals with solutions obtained 
as limits of solutions}
to the $\eps$-problem as $\eps\seto0$.
We also observe that there is no reason for the function $\rhos$ 
mentioned in the statement
to be a constant,
since $f$ is not required to be convex.

\medskip

The rest of the paper is \organiz ed as follows:
in the next section, we 
prove both Proposition~\ref{Globalexistence} and Theorem~\ref{Convergence};
the proof of Theorem~\ref{Omegalimit}
 {will be given}
in the last section.


\section{Global solutions}
\label{ConvProof}
\setcounter{equation}{0}

We first prove the existence of a global solution to the $\eps$-problem.
The  {major} part of the section is devoted to the proof of Theorem~\ref{Convergence}
and the subsequent existence of a global solution to the limit problem.

\step Proof of Proposition~\ref{Globalexistence}

We {imitate, with minor changes, the proof of Thm~2.1 in \cite{CGPS3},}
where the final time $T$ was fixed once and for all.
{Let $\eps$ be fixed and, for notational conciseness, let the dependence on $\eps$ be omitted.}
The main tool used in \cite{CGPS3} was an approximation {procedure} using a time delay $\tau=T/N$,
{for} $N$ a positive integer.
Approximating $\tau$-problems were constructed and solved step by step
in the time intervals $I_n:=[0,n\tau]$, $n=1,\dots,N$.
It turned out that the resulting unique solution $(\mu^\tau,\rho^\tau)$ 
coincided a posteriori with the one obtained
by glueing  {together} solutions on the time steps 
$[(n-1)\tau,n\tau]$, $n=1,\dots,N$.

The  {necessary} slight modification is the following.
Take, e.g., $\tau=1/N$, and solve the same problems as before step by step,
now {for} every $n\geq1$.
This provides a global solution $(\mu^\tau,\rho^\tau)$ 
to the approximating $\tau$-problem.
Then, for every fixed~$T>0$, the argument of \cite{CGPS3} applies,
and a solution for \Pbleps\ on $[0,T]$ is constructed as the limit of the approximating solutions
as $\tau\seto0$, at least for a subsequence.
This holds, in particular, for $T=1,2,\dots\,$.  {Therefore,}
 {there is} a~subsequence $\tau_{1,n}\seto0$ such that
the restriction of $(\mu^\tau,\rho^\tau)$ to $[0,1]$
converges to a solution to problem \Pbleps\ with $T=1$.
We 
 {denote this solution by}
 $(\mu_1,\rho_1)$.
Now, take the restriction of $(\mu^\tau,\rho^\tau)$ to $[0,2]$
with $\tau=\tau_{1,n}$.
Then, for the same reason, {there is}
a~subsequence $\{\tau_{2,n}\}$ of $\{\tau_{1,n}\}$ such that
the restriction we are considering
converges to a solution $(\mu_2,\rho_2)$ to problem \Pbleps\ with $T=2$.
However, as $\{\tau_{2,n}\}$ is a subsequence of~$\{\tau_{1,n}\}$,
the restriction of $(\mu_2,\rho_2)$ to $[0,1]$ coincides with $(\mu_1,\rho_1)$.
 {Proceeding inductively in this way}, 
and then using a diagonal procedure, leads to a global solution to problem \Pbleps.
\QED

\medskip

\step {Preliminaries to the proof of Theorem~\ref{Convergence}}

{We begin by listing some of the tools we shall use.} First of all,  the \wk\ continuous embedding, with the related Sobolev inequality, holds
in our $3$-dimensional case:
\Bsist
  && \Wx{1,p} \subset \Lx q
  \aand
  \norma v_{\Lx q} \leq C_p \norma v_{\Wx{1,p}}
  \quad \hbox{for every $v\in\Wx{1,p}$,}\quad
  \label{gensobolev}
  \\
  && \quad \hbox{provided that} \quad
  1 \leq p <3
  \aand
  1 \leq q \leq p^* := \frac {3p} {3-p} \,,
  \label{sobexp}
\Esist
with the constant $C_p$ in \eqref{gensobolev} depending only on $\Omega$ and~$p$; moreover,
\Beq
  \hbox{the embedding} \quad
  \Wx{1,p} \subset \Lx q
  \quad \hbox{is compact if} \quad
  1 \leq q < p^* .
  \label{compsobolev}
\Eeq
In particular, $V\subset\Lx q$ for $1\leq q\leq 6$, and 
\Beq
  \norma v_{\Lx q} \leq C \normaV v
  \quad \hbox{for every $v\in V$ and $1\leq q \leq 6$,}
  \label{sobolev}
\Eeq
where $C$ depends 
 {only on~$\Omega$ and 
the embedding $V\subset\Lx q$ is compact if $q<6$.}
Furthermore
(see, e.g., \cite[formula~(3.2), p.~8]{DiBen}),
 {we have}
the continuous embedding
\[
 \L\infty H \cap \L2 V \subset \LQ{10/3}
\]
and the related inequality
\Bsist
  && \norma v_{\LQ{10/3}} \leq C_T \norma v_{\L\infty H\cap\L2V}
  \non
  \\
  && \quad \hbox{for every $v\in\L\infty H\cap\L2V$},
  \label{diben}
\Esist
where $C_T$ depends on $\Omega$ and $T$.
In our proof, we shall make use also of
the \wk\ \holder\ inequality, mainly in the form
\Bsist
 &&\norma {v_1\cdots v_n}_{\L p{\Lx q}}
  \leq \prod_{i=1}^n \norma{v_i}_{\L{p_i}{\Lx{q_i}}}\quad\textrm{for}\;\; v_i \in\L{p_i}{\Lx{q_i}}, \;i=1,\dots,n,
  \non\\
   &&\qquad\hbox{provided that} \quad p,q,p_i,q_i \in [1,+\infty] ,
  \quad \frac 1p = \sum_{i=1}^n \frac 1 {p_i}
  \aand
  \frac 1q = \sum_{i=1}^n \frac 1 {q_i} \,.
  \non
\Esist

\Brem
{To avoid a cumbersome notation, the lowercase letter $c$ stands for different constants, each of which may depend on one or another of the data involved in our current statement and on the coefficient~$\delta$, but never depends either on~$\eps$ or on the final time~$T$; consequently, the relative estimates continue to hold when we discuss both the system's asymptotic limit as $\varepsilon$ tends to 0 and its \loti\ \bhv. Moreover,
a~notation like~$c_\sigma$ signals that that constant has an additional dependence on the parameter~$\sigma$.}
Hence, the meaning of $c$ and $c_\sigma$ may
change from line to line, and even in the same chain of inequalities.
On the contrary, {we use the uppercase letter $C$ for precise constants we are going to refer to after their introduction, such as $C_p$ in \eqref{gensobolev} or $C_T$ in \eqref{sobolev}.}
Finally, in order to {lighten our} notation,
we do not write the subscript $\eps$
in performing our a~priori estimates until each estimate is completely proved; {the same we do}
 for the auxiliary function
\Beq
  \ueps := \eps\mueps + 2\mueps\rhoeps.
\label{defueps} 
\Eeq
\Erem

Next, we prove that $\rhoeps$ is bounded away from zero 
uniformly with respect to~$\eps$.
Such a result is essentially known
from the proof of \cite[Thm~2.3]{CGPS3}, 
among other properties there established for a fixed~$\eps$.
Nevertheless, we prefer to repeat the proof here,
in order to  {make sure that the constructed 
lower bound is in fact independent of~$\eps$},
and that just the additional assumption \eqref{hpinfrhoz} is used. 

\Blem
\label{Farfromzero}
There exists some $\rmin\in(0,1)$ such that
$\rhoeps\geq\rmin$ a.e.\ for every $\eps\in(0,1)$.
\Elem

\par\noindent{\bf Proof.}
We set for convenience $\rhomin:=\inf\rhoz$ and $M:=\sup_{r\in(0,1)}|f_2'(r)|$
and recall that $\rhomin > 0$ by~\eqref{hpinfrhoz}.
Thus, owing to~\eqref{hpfp},
we can choose $\rmin\in(0,\rhomin { ]}$ such that $f_1'(\rmin)\leq-M$.
Then, we test \eqref{seconda} by $-\rhomrmeps$ and integrate over $\Omega\times(0,t)$
where $t\in(0,T)$ is arbitrary.
By omitting the subscript $\eps$ for simplicity, we have
\Bsist
  && \frac \delta 2 \iO |\rhomrm(t)|^2
  + \intQt |\nabla\rhomrm|^2
  - \intQt \bigl( f_1'(\rho) - f_1'(\rmin) \bigl) \rhomrm
  \non
  \\
  && = \frac \delta 2 \iO |\rhomrm(0)|^2
  - \intQt \mu \rhomrm 
  + \intQt \bigl( f_1'(\rmin) - f_2'(\rho) \bigl) \rhomrm .
  \non
\Esist
Every term on the \lhs\ is nonnegative; {in the \rhs,
the first term vanishes, because $\rhoz\geq\rmin$, and the other
 two are nonpositive, because} $\mu\geq0$ and $f_1'(\rmin)-f_2'(\rho)\leq f_1'(\rmin)+M\leq0$.
Hence, $\rhomrm=0$,
{and the assertion is proved.}
\QED

\step {Proof of Theorem~\ref{Convergence}}

Our proof will proceed as follows.
For a fixed finite final time~$T$, we shall perform a number of a~priori estimates
and use \wk\ compactness results  
{to prove that, as $\eps$ tends to 0,}
the solution $(\mueps,\rhoeps)$ to the $\eps$-problem \Pbleps\ 
we are considering converges
to a solution $(\mu,\rho)$ to problem \Pbl, at least for a subsequence $\eps_n\seto0$; in particular, this holds for $T=1,2,\dots\,$.
{Having established this result,}
 we shall be able to argue as in the proof of Proposition~\ref{Globalexistence}.
Indeed, a~diagonal procedure provides a subsequence $\eps_n\seto0$
such that $(\mueps,\rhoeps)$ with $\eps=\eps_n$ converges
to a global solution $(\mu,\rho)$ to problem \Pbl\
defined in the whole of~$[0,+\infty)$.
Therefore, just the case of a fixed final time $T$ has to be considered.



\step First a priori estimate

We test~\eqref{primaeps} (e.g., the equation {within square brackets})
by~$\mueps$, and integrate over $\Omega\times(0,t)$,
for an arbitrary $t\in(0,T)$.
We obtain
\Beq
  \intQt \dt\!\left(  \frac\eps 2 \, \mu^2 + \rho\mu^2 \right)
  + \intQt |\nabla\mu|^2 = 0,
  \quad 
  \non
\Eeq  
whence
\Beq
  \frac \eps 2 \, \iO |\mu(t)|^2
  + \iO (\rho\mu^2)(t)
  + \intQt |\nabla\mu|^2 
   = \frac \eps 2 \normaH\muz^2
  + \norma{\rhoz\muz^2}_{\Luno}
  \leq c .
  \non
\Eeq
Since $\rho\mu^2\geq\rmin\mu^2$ thanks to Lemma~\ref{Farfromzero},
we immediately deduce that
\Beq
  \norma\mueps_{\L\infty H}
  + \norma{\nabla\mueps}_{\L2H}
  \leq c .
  \label{primastima}
\Eeq

\step Second a priori estimate

We test \eqref{secondaeps} by $\dt\rhoeps$,
and use the second {of} \eqref{primaeps}
in order to compute the \rhs\ we get; we also recall~\eqref{defueps}.
For $t\in(0,T)$, we obtain:
\Bsist
  && \delta \intQt |\dt\rho|^2
  + \frac 12 \iO |\nabla\rho(t)|^2
  - \frac 12 \iO |\nabla\rhoz|^2
  + \iO f(\rho(t))
  - \iO f(\rhoz)
  \non
  \\
  && = \intQt \mu \dt\rho
  = \intQt \dt \bigl( \eps\mu + 2\mu\rho \bigr)
  - \intQt \Delta\mu
  \non
  \\
  && = \intQt \dt u
  = \iO u(t) - \iO \bigl( \eps\muz + 2\rhoz\muz \bigr)
  \leq 3 \iO \mu(t) + c .
  \non
\Esist
Since \eqref{primastima} holds and $f$ is bounded from below, 
we easily infer that
\Beq
  \norma{\dt\rhoeps}_{\L2H}
  + \norma{\nabla\rhoeps}_{\L\infty H}
  + \norma{f(\rhoeps)}_{\L\infty\Luno}
  \leq c .
  \label{secondalong}
\Eeq
{Moreover, because}
 $0<\rhoeps<1$ \aeQ\ for every $\eps\in(0,1)$,
we conclude that
\Beq
  \norma\rhoeps_{\H1H\cap\L\infty V} \leq c_T.
  \label{secondastima}
\Eeq

\step Third a priori estimate

{Taking into account} 
\accorpa{primastima}{secondastima}
and the boundedness of~$f'_2$,
we see that \eqref{secondaeps} yields
\Beq
  \norma{-\Delta\rhoeps + f_1'(\rhoeps)}_{\L2H}
  = \norma{\mueps - \dt\rhoeps - f_2'(\rhoeps)}_{\L2H}
  \leq c_T \,.
  \non
\Eeq
By a standard argument
(test formally by~$f_1'(\rhoeps)$, for instance)
and~elliptic regularity, we conclude~that 
\Beq
  \norma{f_1'(\rhoeps)}_{\L2H} + \norma\rhoeps_{\L2W} \leq c_T \,.
  \label{terzastima}
\Eeq

\step First conclusions

The above estimates allow us to use
standard weak and weak star compactness results.
Thus, a~triplet $(\mu,\rho,\phi)$ exists such~that
\Bsist
  && \mueps \to \mu
  \quad \hbox{weakly star in $\L\infty H\cap\L2V$}
  \label{convmu}
  \\
  && \rhoeps \to \rho 
  \quad \hbox{weakly star in $\H1H\cap\L\infty V\cap\L2W$}
  \label{convrho}
  \\
  && f_1'(\rhoeps) \to \phi
  \quad \hbox{weakly in $\L2H$},
  \label{convfp}
\Esist
at least for a subsequence $\eps_n\seto0$.\footnote{{Incidentally, we anticipate that} the convergence results 
{stated below will hold only}
for suitable subsequences.
Nevertheless, we will not mention such a detail.}
We note  that $\mu\geq0$ and that $\rho\geq\rmin$ \aeQ\
(the former {inequality holds because} $\mueps\geq0$ for every~$\eps$,
the latter by Lemma~\ref{Farfromzero}).
Moreover, by \eqref{convrho} and the compact embedding $V\subset\Lx p$ for $p<6$,
we~infer that
\Beq
  \rhoeps \to \rho
  \quad \hbox{strongly in $\C0{\Lx p}$}
  \quad \hbox{for every $p<6$},
  \label{strongrho}
\Eeq
thanks to \cite[Sect.~8, Cor.~4]{Simon}.
{Hence}, 
$f_2'(\rhoeps)$ converges to $f_2'(\rho)$ in a suitable topology, 
for instance, strongly in $\L2H$, since $f_2'$ is Lipschitz continuous.
In particular, we deduce that
\Beq
  \delta \dt\rho - \Delta\rho + \phi + f_2'(\rho) = \mu 
  \quad \aeQ .
  \non
\Eeq
Furthermore, 
{invoking}
both \eqref{strongrho} and~\eqref{convfp}, 
and using a standard monotonicity technique
(see, e.g., \cite[Lemma~1.3, p.~42]{Barbu}),
we conclude~that
\Beq
  0 < \rho < 1
  \aand
  \phi = f_1'(\rho)
  \quad \aeQ .
  \non
\Eeq
Finally, \eqref{strongrho} implies that $\rhoeps(0)$ converges to~$\rho(0)$ strongly in~$H$,
whence $\rho(0)=\rhoz$.

\medskip

{In summary,} {so far}
we have proved \accorpa{regmu}{regfprho}, \eqref{seconda}, and 
{the second condition in~\eqref{cauchy}. It}
remains {for us} to show \eqref{regprod}, \eqref{prima}, and the first
{condition in~\eqref{cauchy}.
For this purpose, further arguments are needed.}

\step Fourth a priori estimate

We recall~\eqref{defueps}, and we notice that the second \eqref{primaeps} reads:
\Beq
  \dt\ueps = \mueps \dt\rhoeps + \Delta\mueps \,.
  \label{reprimaeps}
\Eeq
Moreover, $\mueps$ satisfies homogeneous Neumann boundary conditions, since it is $W$-valued.
Therefore, we have
\Beq
  \intQ \dt\ueps \, v
  = \intQ \mueps \dt\rhoeps \, v
  - \intQ \nabla\mueps \cdot \nabla v
  \quad \hbox{for every $v\in\L2V$}.
  \label{varprimaeps}
\Eeq
From \eqref{varprimaeps}, we derive an estimate for $\dt\ueps$
as a $\Vp$-valued function,
in the framework of the Hilbert triplet $(V,H,\Vp)$.
We treat the integrals on the \rhs\
{individually}
{(for a while, we} omit the subscript $\eps$ in order to simplify the notation).

{Assume that}
 $v\in\L5V$.
Then the \holder\ inequality and the Sobolev inequality \eqref{sobolev} with $q=5$ yield:
\Bsist
  && \Bigl|
    \intQ \mu \dt\rho \, v
  \Bigr|
  \leq \norma\mu_{\LQ{10/3}} \, \norma{\dt\rho}_{\LQ2} \, \norma v_{\LQ5} 
  \non
  \\
  && \leq { c} \norma\mu_{\LQ{10/3}} \, \norma{\dt\rho}_{\LQ2} \, \norma v_{\L5V} \,.
  \non
\Esist
On the other hand, inequality \eqref{diben} holds.
Therefore, {on taking into account} \eqref{primastima} and~\eqref{secondastima},
we conclude that
\Beq
  \Bigl|
    \intQ \mu \dt\rho \, v
  \Bigr|
  \leq c_T \norma v_{\L5V}
  \quad \hbox{for every $v\in\L5V$}.
  \label{perdtu}
\Eeq

{Next, we}
consider the second term on the \rhs\ of~\eqref{varprimaeps}.
By assuming 
{again that $v\in\L5V$, and invoking~\eqref{primastima} once more, we immediately find that}
\Beq
  \Bigl|
    \intQ \nabla\mu \cdot \nabla v
  \Bigr|
  \leq \norma\mu_{\L2V} \, \norma v_{\L2V}
  \leq c_T \norma v_{\L5V} \,.
  \non
\Eeq
Combining 
{this estimate with \eqref{perdtu} and \eqref{varprimaeps},
we obtain that}
\Beq
  | \< \dt u , v > |
  \leq c_T \norma v_{\L5V}
  \quad \hbox{for every $v\in\L5V$}.
  \non
\Eeq
{In other words,} {we have that}
\Beq
  \norma{\dt\ueps}_{\L{5/4}\Vp} \leq c_T \,.
  \label{quartastima}
\Eeq

\step Consequence

From {(i)} the strong convergence \eqref{strongrho} with $p=4$, {(ii)} 
the weak convergence $\mueps\to\mu$ in $\L2\Lq$
implied by~\eqref{convmu}, 
and {(iii)}  the Sobolev inequality \eqref{sobolev} with $q=4$, 
we infer that
\Beq
  \mueps \rhoeps \to \mu \rho
  \quad \hbox{weakly in $\L2H$},
  \  \, \hbox{whence} \  \,
  \ueps \to 2\mu\rho
  \ \, \hbox{weakly in $\L2H$},
  \label{convu}
\Eeq
since $\eps\mueps\to0$ strongly in $\L2V$, by~\eqref{primastima}.
Hence, accounting for \eqref{quartastima}, we conclude~that
\Beq
\dt\ueps \to 2\dt(\mu\rho)
  \ \, \hbox{weakly in $\L{5/4}\Vp$}, \  \, \hbox{whence} \  \,\ueps \to 2\mu\rho
  \ \, \hbox{weakly in $\W{1,5/4}\Vp$},  \label{moreu}
 \Eeq
  %
so that \eqref{regprod} holds with $p=5/4$.
Moreover, \eqref{moreu}~also implies that
$\ueps$ converges to $2\mu\rho$ weakly in $\C0\Vp$; in particular, $\ueps(0)\to(2\mu\rho)(0)$ weakly in $\Vp$.
On the other hand,
$\ueps(0)=\eps\muz+2\muz\rhoz$ converges to $2\muz\rhoz$, e.g., strongly in~$H$.
Thus, the Cauchy condition for $\mu\rho$ in \eqref{cauchy}  follows.

\medskip

In order to prove~\eqref{prima},
one can try to let $\eps\seto0$ in \eqref{varprimaeps} first, then to get rid of time integration. But, a trouble arises in dealing with the first term on the right-hand side,
{since, for the moment being,
 both $\mueps$ and $\dt\rhoeps$ are just weakly convergent.}
Hence, we have to prepare a new tool.

\step Fifth a priori estimate

We want to find a bound for $\nabla\ueps$, i.e.,
for the partial derivatives $D_i\ueps$, $i=1,2,3$.
As usual, we omit the subscript $\eps$ for a while.
We have:
\Beq
  |D_i u|
  = |\eps D_i\mu + \rho D_i\mu + \mu D_i\rho|
  \leq 2 |D_i\mu| + \mu |D_i\rho|.
  \non
\Eeq
Now, on
{taking \eqref{primastima} into account,}
we see that $D_i\mu$ is bounded in $\L2H$,
{while}
$\mu$~is bounded in $\L2{\Lx6}$ thanks to the Sobolev inequality~\eqref{sobolev}.
On the other hand, 
{\eqref{secondastima}}
provides a bound for $D_i\rho$ in $\L\infty H$.
Hence, using \holder\ inequality,
we see that the product $\mu D_i\rho$ is bounded in $\L2{\Lx{3/2}}$.
Therefore, we conclude that
\Beq
  \norma\ueps_{\L2{\Wx{1,3/2}}} \leq c_T \,.
  \label{quintastima}
\Eeq

\step Consequence

As \eqref{convu} holds, from \eqref{quintastima} we infer that
\Beq
  \ueps \to 2\mu\rho
  \quad \hbox{weakly in $\L2{\Wx{1,3/2}}$}.
  \label{betteru}
\Eeq
Now, we observe that the embedding $\Wx{1,3/2}\subset\Lx q$ is compact for every $q<3$,
by~\eqref{compsobolev}.
On the other hand, \eqref{moreu} holds.
By using the Aubin-Lions lemma 
(see, e.g., \cite[Thm.~5.1, p.~58]{Lions}),
we deduce the strong convergence
\Beq
  \ueps \to 2\mu\rho
  \quad \hbox{strongly in $\L2{\Lx q}$}
  \quad \hbox{for every $q<3$}.
  \label{strongu}
\Eeq
We stress that, in particular, $\ueps\to2\mu\rho$ strongly in $\L2H$.

\Blem
\label{Strongmu}
The strong convergence $\mueps\to\mu$ holds in $\L2H$.
\Elem

\par\noindent{\bf Proof.}
We set $u:=2\mu\rho$ and argue \aeQ\ for a while.
Thanks to Lemma~\ref{Farfromzero}, we have
\Beq
  |\mueps - \mu|
  = \Bigl|
    \frac \ueps {\eps+2\rhoeps} - \frac u {2\rho}
  \Bigr|
  = \Bigl|
    \frac {2\rho\ueps - \eps u - 2 \rhoeps u} {2\rho (\eps + 2\rhoeps)}
  \Bigr|
  \leq \frac {\eps |u| + 2 |\rho\ueps - \rhoeps u|} {4\rmin^2} \,.
  \non
\Eeq
On the other hand, we have
\Beq
  |\rho\ueps - \rhoeps u|
  \leq |\rho| \, |\ueps - u| + |u| \, |\rho - \rhoeps|
  \leq |\ueps - u| + 2\mu \, |\rho - \rhoeps| .
  \non
\Eeq
By combining {these inequalities}, we deduce that
\Bsist
  \hskip-1.5cm&& \norma{\mueps-\mu}_{\L2H} \non\\
  \hskip-1.5cm&& \leq c \left(
    \eps \norma u_{\L2H}
    + \norma{\ueps-u}_{\L2H}
    { + \norma{\mu}_{{{L^2(0,T; L^4(\Omega))}}}
      \norma{\rho-\rhoeps}_{C^0([0,T]; L^4(\Omega){)}}}
  \right) .
  \label{perstrongmu}
\Esist
The first two terms on the \rhs\ tend to zero as $\eps\seto0$,
by 
\eqref{strongu}{; as to
the last term,}
it suffices to recall \eqref{regmu}, \eqref{strongrho}, and \eqref{sobolev} for $q=4$.
\QED

\step Conclusion

The strong convergence {guaranteed} by Lemma~\ref{Strongmu}, {together with}
the weak convergence $\dt\rhoeps\to\dt\rho$ in $\L2H$
given by \eqref{convrho}, imply that
\Beq
  \mueps \, \dt\rhoeps \to \mu \, \dt\rho
  \quad \hbox{weakly in $\LQ1$}.
\non
\Eeq
On the other hand, \eqref{moreu} and \eqref{convmu} hold.
Hence, by letting $\eps\seto0$ in \eqref{varprimaeps}, we easily obtain
that
\Bsist
  && 2 \ioT \< \dt(\mu\rho)(t) , z(t) > \, dt
  = \intQ \mu \dt\rho \, z
  - \intQ \nabla\mu \cdot \nabla z
  \label{varprima}
\Esist
for every $z\in\L5V\cap\LQ\infty$.
Now, take any $v\in V\cap\Lx\infty$ and any $\zeta\in L^\infty(0,T)$.
Then the function $z:t\mapsto\zeta(t)v$ is admissible in~\eqref{varprima},
and a standard argument yields
\Beq
  2 \< \dt (\mu\rho)(t) , v >
  + \iO \nabla\mu(t) \cdot \nabla v
  = \iO \mu(t) \, \dt\rho(t) \, v,
  \non
\Eeq
for every $v\in V\cap\Lx\infty$ and \aat.
Now, we note that, \aat, each term in the above equation
defines an element of~$\Vp$.
This is clear as far the \lhs\ is concerned,
since $\dt(\mu\rho)$ is $\Vp$-valued and $\mu$ is $V$-valued.
For the remaining term, we recall Remark~\ref{Commentoprima}.
On the other hand, $V\cap\Lx\infty$ is dense in~$V$.
Therefore, \eqref{prima}~follows,
and the proof is complete.
\QED


\section{\Loti\ \bhv}
\label{OmegaProof}
\setcounter{equation}{0}

This section is devoted to {proving} Theorem~\ref{Omegalimit}.
We first derive some a~priori estimates, then  we prove the continuity property announced in the statement; finally, we \characteriz e the $\omega$-limit.

\step A priori estimates

We recall Lemma~\ref{Farfromzero} and the a priori estimates 
\eqref{primastima} and~\eqref{secondalong},
which involve constants that do not depend on the final time.
Hence, we immediately obtain that
\Beq
  \rho(t) \geq \rmin 
  \aand
  \normaH{\mu(t)} 
  + \normaH{\nabla\rho(t)}
  \leq c 
  \quad \Aat .
  \label{stimelongtime}
\Eeq
Recalling that $0<\rho<1$, we see, in particular, 
that $(\mu,\rho)$ is a bounded $(H\times V)$-valued function,
as stated.
For the same reason, we also deduce that
\Beq
  { \norma{(\mu\rho)(t)}_H} \leq c
  \quad \Aat .
  \label{ubdd}
\Eeq
Moreover, 
{the estimates \eqref{primastima} and \eqref{secondalong}
also yield the bounds}
\Beq
  \norma{\nabla\mu}_{\L2H}
  + \norma{\dt\rho}_{\L2H}
  \leq c
  \quad \hbox{for every $T>0$},
  \non
\Eeq
and we conclude that
\Beq
  \intQi |\nabla\mu|^2 < +\infty
  \aand
  \intQi |\dt\rho|^2 < +\infty .
  \label{spacetimeindep}
\Eeq

\step Strong and weak continuity

As far as $\rho$ is concerned,
we have $\rho\in\H1H\cap\L2W$ for every $T<+\infty$ by~\eqref{convrho}.
Since the embedding 
\Beq
\H1H\cap\L2W\subset\C0V 
\non
\Eeq
holds,
we immediately deduce that $\rho$ is a strongly continuous $V$-valued function.
The weak continuity of $\mu$ is less obvious: we prove it by using the following \wk\ tool, whose proof {is left as an exercise to the reader}.

\Bprop 
\label{Tool}
Let~$\calZ$ be a Hausdorff topological space,
and let $Z$ be a reflexive Banach space such that $Z\subset\calZ$,
where the embedding is continuous with respect to the weak topology of~$Z$.
Assume that $z:[0,T]\to\calZ$ is continuous and that
$z(t)\in Z$ and $\norma{z(t)}_Z\leq M$ \aat\ for some constant~$M$.
Then $z$~is $Z$-valued, i.e., $z(t)$ {belongs to} $Z$ for every $t\in[0,T]$,
and {is} continuous with respect to the weak topology of~$Z$.
Moreover, $\norma{z(t)}_Z\leq M$ for every $t\in[0,T]$.
\Eprop

In our case, we argue on any fixed finite time interval $[0,T]$ 
and apply {Proposition~\ref{Tool}} twice,
first with $\calZ = { V^*}$,
with either the weak or the strong topology,
then with $\calZ=\Luno$, endowed with the weak topology.
We set: 
\Beq
  u := 2\mu\rho,
  \label{defu}
\Eeq 
in order to agree with \eqref{defueps}, 
and we recall that $u\in\W{1,5/4}\Vp$, by~\eqref{moreu}; in particular, $u\in\C0\Vp$.
On the other hand, we have proved~\eqref{ubdd}.
We conclude that $u(t)\in H$  for every $t\in[0,T]$
and that $u$ is continuous with respect to the weak topology of~$H$.
Besides, $\rho$~is strongly continuous as an $H$-valued function,
and the first 
{condition in}
\eqref{stimelongtime} holds.
Hence, the same is true for~$1/\rho$,
and we infer that $\mu=u/(2\rho)$ is a weakly continuous $\Luno$-valued function.
Now, we recall the estimate of $\mu$ given by~\eqref{stimelongtime}
and conclude that $\mu$ is weakly continuous as an $H$-valued function.

\step Conclusion

It remains {for us} to show that every element of the $\omega$-limit
is a steady state.
{To this end,}
we pick any $(\muo,\rhoo)\in\omega(\mu,\rho)$ 
and consider a corresponding sequence $t_n\neto+\infty$, as
given by definition~\eqref{defomegalim}.
We~set:
\Beq
  \mun(t) := \mu(t_n+t) , \quad
  \rhon(t) := \rho(t_n+t) ,
  \aand
  \un(t) := u(t+t_n),
  \quad \hbox{for $t\geq0$,}
  \label{munrhon}
\Eeq
and study the sequence $\graffe{(\mun,\rhon)}$ 
on a fixed finite time interval~$[0,T]$
by using $\un$ as well.
From~\eqref{stimelongtime}, \eqref{spacetimeindep},
and weak star compactness, we immediately deduce that
\Bsist
  && \mun \to \mui
  \quad \hbox{weakly star in $\L\infty H$},
  \quad  \rhon \to \rhoi
  \quad \hbox{weakly star in $\L\infty V$,}
  \non
  \\
  && |\nabla\mun| \to 0
  \aand 
  \dt\rhon \to 0
  \quad \hbox{strongly in $\L2H$,}
  \non
\Esist
at least for a subsequence.
It follows that $\mui$ is {space-
and $\rhoi$ time-}independent.
Thus, we can write $\rhoi(t)=\rhos$ \aat\
for some $\rhos\in V$.
On the other hand,
we can reproduce 
{the}
estimates \eqref{terzastima}, \eqref{quartastima}, and \eqref{quintastima},
on the time interval $[t_n,t_n+T]$
instead of~$[0,T]$.
We obtain:  
\Beq
  \norma{f_1'(\rho)}_{L^2(t_n,t_n+T;H)}
  + \norma\rho_{L^2(t_n,t_n+T;W)}
  + \norma\un_{W^{1,5/4}(t_n,t_n+T;\Vp)\cap L^2(t_n,t_n+T;\Wx{1,3/2})}
  \leq c_T,
  \non
\Eeq
{where}
 $c_T$ does not depend on~$n$.
This means that
\Beq
  \norma{f_1'(\rhon)}_{\L2H}
  + \norma\rhon_{\L2W}
  + \norma\un_{\W{1,5/4}\Vp\cap\L2{\Wx{1,3/2}}}
  \leq c_T \,.
  \label{perforte}
\Eeq
Thus, $\rhoi\in\L2W$, i.e., $\rhos\in W$.
Moreover, 
to derive a strong convergence for $\rhon$ in~$\C0H$, we can argue as in the previous section.
This allows us to ensure that 
$f_2'(\rhon)$ converges to $f_2'(\rhoi)$ strongly in $\L2H$
and that the weak limit of $f_1'(\rhon)$ in $\L2H$, 
given by weak compactness, is~$f_1'(\rhoi)$.
All this yields that
\Beq
  0 < \rhos < 1
  \aand
  - \Delta\rhos + f'(\rhos) = \mui
  \quad \aeQ ,
  \non
\Eeq
and we deduce that $\mui$ is even time-independent.
Thus, $\mui(x,t)=\mus$ \aaQ\ for some constant~$\mus$.
Furthermore, $\mus$~is nonnegative, since $\mun\geq0$ for every~$n$.
This concludes the proof that $(\mus,\rhos)$ is a steady state.

{Lastly, we} show that $(\mus,\rhos)$ coincides with~$(\muo,\rhoo)$.
Because $\rhon\to\rho$ strongly in $\C0H$, we see that
$\rhon(0)$ converges to~$\rhoi(0)=\rhos$ strongly in~$H$.
On the other hand,  by assumption $\rhon(0)=\rho(t_n)$ converges to $\rhoo$
weakly in $V$.
Hence, $\rhos=\rhoo$.
The 
{corresponding}
argument for $\mus$ and $\muo$ 
{is a bit more involved.}
We remind that the embedding $\Wx{1,3/2}\subset H$ is compact.
Hence, from \eqref{perforte} and the Aubin-Lions lemma, 
{we conclude that there is some  $\ui$ such~that}
\Bsist
  && \un \to \ui 
  \quad \hbox{weakly in $\W{1,5/4}\Vp\cap\L2{\Wx{1,3/2}}$,}
  \label{weakun}
  \\
  && \hbox{whence} \quad
  \un \to \ui 
  \quad \hbox{strongly in $\L2H$}.
  \label{strongun}
\Esist
On the other hand, the strong convergence $\rhon\to\rhoi$ in $\C0H$
and the uniform inequality $\rhon\geq\rmin$
imply the strong convergence $1/\rhon\to1/\rhoi$ in $\C0H$.
We infer that
\Beq
  \mun = \frac \un {2\rhon} \to \frac \ui {2\rhoi}
  \quad \hbox{strongly in $\L2\Luno$}.
  \non
\Eeq
Since $\mun\to\mui$ weakly star in $\L\infty H$,
we conclude that $\ui/(2\rhoi)=\mui$, i.e., that
$\ui(t)=2\mus\rhoo$ \aat.
Next, the first weak convergence \eqref{weakun} also implies
weak convergence in $\C0\Vp$; in particular, $\un(0)$~converges to $\ui(0)=2\mus\rhoo$ weakly in~$\Vp$.
On the other hand, by assumption
$\mu(t_n)\to\muo$ weakly in~$H$ 
{and, due} to the already mentioned strong convergence
$\rhon\to\rhoi$, $\rho(t_n)\to\rhoo$ strongly in~$H$.
We infer that $\un(0)=2\mu(t_n)\rho(t_n)$ converges to $2\muo\rhoo$
weakly in~$\Luno$.
By comparison, we conclude that
$2\mus\rhoo=2\muo\rhoo$,
i.e., that $\mus=\muo$.
\QED



\vspace{3truemm}

\Begin{thebibliography}{99}

\bibitem{Barbu}
V. Barbu,
``Nonlinear semigroups and differential equations in Banach spaces'',
Noord\-hoff,
Leyden,
1976.

\bibitem{CN} B.D. Coleman, W. Noll, 
The thermodynamics of elastic materials with heat conduction and viscosity, 
{\it Arch. Rational Mech. Anal.} {\bf 13} (1963) 167-178.

\bibitem{CGPS} 
P. Colli, G. Gilardi, P. Podio-Guidugli, J. Sprekels,
Existence and uniqueness of a global-in-time solution
to a phase segregation problem of the Allen-Cahn type, 
{\it Math. Models Methods Appl. Sci.} {\bf 20} (2010)
519-541.

\bibitem{CGPS2} 
P. Colli, G. Gilardi, P. Podio-Guidugli, J. Sprekels,
A temperature-dependent phase segregation problem 
of the Allen-Cahn type, {\it Adv. Math. Sci. Appl.} 
{\bf 20} (2010) 219-234.

\bibitem{CGPS3}
P. Colli, G. Gilardi, P. Podio-Guidugli, J. Sprekels,
\Wepo ness and \loti\ \bhv\ for a nonstandard 
viscous Cahn-Hilliard system, WIAS Preprint No. 1602, 
Berlin 2011; to appear in {\it SIAM J. Appl. Math.}

\bibitem{DiBen}
E. DiBenedetto,
``Degenerate Parabolic Equations'',
Springer-Verlag,
New York,
1993.

\bibitem{Fremond}
M. Fr\'emond,
``Non-smooth Thermomechanics'',
Springer-Verlag, Berlin, 2002.

\bibitem{FG} 
E. Fried, M.E. Gurtin, 
Continuum theory of thermally induced phase transitions based on an order 
parameter, {\it Phys. D} {\bf 68} (1993) 326-343.

\bibitem{Gurtin} 
M. Gurtin, Generalized Ginzburg-Landau and
Cahn-Hilliard equations based on a microforce balance,
{\it Phys.~D\/} {\bf 92} (1996) 178-192.

\bibitem{Lions}
J.L. Lions,
``Quelques m\'ethodes de r\'esolution des probl\`emes aux limites non
lin\'eaires'',
Dunod Gauthier--Villars,
Paris,
1969.

\bibitem{Podio}
P. Podio-Guidugli, 
Models of phase segregation and diffusion of atomic species on a lattice,
{\it Ric. Mat.} {\bf 55} (2006) 105-118.

\bibitem{Simon}
J. Simon,
{Compact sets in the space $L^p(0,T; B)$},
{\it Ann. Mat. Pura Appl.} {\bf 146} (1987) 65--96.

\End{thebibliography}

\End{document}

\bye